\documentclass[preprint,12pt]{elsarticle}

\usepackage{amsmath}
\usepackage{amsthm}
\usepackage{amscd}
\usepackage{amssymb}
\usepackage{amsfonts}
\usepackage{hyperref}

\newcommand{\F}{\mathbb{F}}

\newtheorem{defi}{Definition}
\newtheorem{lem}{Lemma}
\newtheorem{thm}{Theorem}
\newtheorem{prop}{Proposition}

\newtheorem{rem}{Remark}
\newtheorem{coro}{Corollary}

\newcommand{\Z}{\mathbb{Z}}

\journal{journal}

\begin{document}

\begin{frontmatter}

\title{The uncertainty principle over finite fields}

\author{Martino Borello}
\address{Universit\'e Paris 8, Laboratoire de G\'eom\'etrie, Analyse et Applications, LAGA,
	Universit\'e Sorbonne Paris Nord, CNRS, UMR 7539, F-93430, Villetaneuse, France}
\author{Patrick Sol\'e}
\address{Aix Marseille University, CNRS, Centrale
Marseille, I2M, Marseille, France}

\begin{abstract}

In this paper we study the uncertainty principle (UP) connecting a function over a finite field and its Mattson-Solomon polynomial,
which is a kind of Fourier transform in positive characteristic. Three versions of the UP over finite fields are studied, in connection with the asymptotic theory of cyclic codes. We first show that no finite field satisfies the \emph{strong} version of UP, introduced recently by Evra, Kowalsky, Lubotzky, 2017. A refinement of the \emph{weak} version is given, by using the asymptotic Plotkin bound. A \emph{naive} version, which is the direct analogue over finite fields of the Donoho-Stark bound over 
the complex numbers, is proved by using the BCH bound.  It is strong enough to show that there exist sequences of cyclic codes of length $n,$ arbitrary rate, and minimum distance $\Omega(n^\alpha)$ for all $0<\alpha<1/2.$
Finally, a connection with Ramsey Theory is pointed out.

\end{abstract}


\end{frontmatter}

\section{Introduction}
The uncertainty principle (UP) is a very famous inequality in Physics \cite{H}, and Signal Processing \cite{D} (see \cite{WW} for a general very recent survey  on the UP). It compares the supports of functions 
and of their complex-valued Fourier transforms.
In a paper of 2017 \cite{EKL}, a connection between UP and the asymptotic performance of cyclic codes was pointed out.
Note that the existence of an asymptotically good family of cyclic codes is a problem open for more than half a century \cite{MPW}. The reference \cite{EKL} is an attempt
to motivate further research into, and eventually solve this very hard problem.
In a recent note \cite{QRS},
a connection with Ramsey Theory and the Szemer\'edi Theorem was derived. 

In the present paper, we replace the classical Discrete Fourier transform (\cite[\S2]{D}) by a vectorial version of the Mattson-Solomon polynomial
(\cite[Ch.8, \S6]{MS}). In contrast with all the results in \cite{WW}, this transform takes its values in a finite field. We study three versions of the UP for this kind of transform. 

The {\em strong} version of the UP over finite fields is defined in \cite{EKL} by analogy with the bound of \cite{T} for the classical Fourier transform. Exploiting the connection with the theory of MDS codes, we show that no finite fields may satisfy the strong UP. 

The {\em weak} version of the UP
is a similar and weakened statement depending on two real parameters $\lambda$ and $\epsilon.$ In \cite{EKL} it is shown that a finite field satisfying this UP enjoys sequences of asymptotically good cyclic codes. Here, we show that, if this version holds over $\F_q,$ then 
$\lambda<\frac{q-1}{q}.$ 

The third version is the straight analogue of the Donoho-Stark bound of \cite{D} and we call this the {\em naive} version. It allows us to construct sequences of cyclic codes with nonzero rate and minimum distance growing
like a power $\alpha$ of the length with $0<\alpha<1/2.$ 

Finally, with similar arguments, we give an alternate proof of 
the results of \cite{QRS}, based on the familiar BCH bound and a generalization based on the Hartmann-Tzeng bound on the minimum distance of cyclic codes (\cite[Th. 4.5.6]{HP}).

\medskip

The material is organized as follows: the next section collects background material; Section 3 is about the strong version; Section 4 contains numerical results related to the weak version; 
Section 5 is dedicated to the naive version of UP; Section 6 deals with the Ramsey Theory 
connection; 
Section 7 concludes the article. An Appendix
building on the naive version shows the existence of cyclic
codes of all rates with minimum distance $\Omega(n^\alpha),$ for all $0<\alpha<1/2.$

\section{Background}

Throughout the paper, $\F_q$ denotes a finite field of cardinality $q$, where $q$ is a prime power.

\subsection{Linear codes and asymptotics}
The {\bf (Hamming) weight}  of $x \in \F_q^n$ is denoted by $w_H(x).$  The minimum nonzero weight $d$ of a linear code is called the {\bf minimum distance}.
A {\bf linear code} is a subspace of $\F_q^n$. Its parameters are written as $[n,k,d]$ where $k$ is its {\bf dimension} as an $\F_q$-vector space.
If $\mathcal{C}_n$ is a sequence of linear codes of parameters $[n, k_n, d_n]$, the {\bf rate} $R$ and {\bf relative distance} $\delta$ are defined as
$$R:=\liminf\limits_{n \rightarrow \infty}\frac{k_n}{n} \ {\rm and} \
\delta:=\liminf\limits_{n \rightarrow \infty}\frac{d_n}{n}.$$
A family of codes is said to be {\bf good} iff it contains a sequence with rate and relative distance such that $R\cdot \delta \neq 0.$
The {\bf binary entropy function} $H(x)$ of the real variable $x$ is defined (see \cite[p.308]{MS}) for $0<x<1,$
$$H(x):=-x\log_2(x)-(1-x)\log_2(1-x).$$

\subsection{Cyclic Codes}
Consider the quotient ring $R(\F_q,n):=\F_q[x]/(x^n-1).$ We will identify each class of $R(\F_q,n)$ with the unique polynomial of degree less than $n$ contained in it. The ring $R(\F_q,n)$ is principal, and we denote by $C(f)$ the ideal with 
generator $f.$ It is well-known that every ideal of $R(\F_q,n)$ has a unique monic generator of minimal degree, and this is a divisor of $x^n-1$. Whenever we will consider and ideal $C(f)$, we will implicitly assume that $f$ is such generator.
The polynomials of $R(\F_q,n)$ are in one-to-one correspondence with the vectors of $\F_q^n,$ by the map

$$ \varphi: f:=(f_0,f_1,\dots, f_{n-1}) \mapsto f(x):=\sum_{i=0}^{n-1} f_i x^i. $$

The (Hamming) weight of a polynomial is the (Hamming) weight of the corresponding vector.
A {\bf cyclic code} is an ideal in $R(\F_q,n)$ or its preimage in $\F_q^n$ via $\varphi$. 

The {\bf zeros} of $C(f)$ are the roots of $f$ in the algebraic closure of $\F_q.$
The dimension of $C(f)$ is $n-\deg(f),$ and $\deg(f)$ equals the number of zeros of $C(f)$ (see for example \cite[Chap. 7]{MS}). The well-known \textbf{BCH-bound} \cite[Chap., Th. 8]{MS} states that if among the zeros of $f$ there exists $\delta-1$ consecutive powers of a primitive $n$-th root of unity and $(m,q)=1$, then the minimum distance of $C(f)$ is at least $\delta$. 

\subsection{Mattson-Solomon polynomial}

Let $\zeta$ be a primitive root of unity of order $n$ in the algebraic closure $\overline{\F_q}$ of $\F_q$. The \textbf{Mattson-Solomon polynomial} (\cite[Chap. 8]{MS}) associated with a vector $f:=(f_0,f_1,\dots, f_{n-1})$ is the following polynomial in $\overline{\F_q}[z]$:
$$\hat{f}(z):=\sum_{i=1}^n F_i z^{n-i},$$
where $F_i:=f(\zeta^i)$ is the evaluation of $f(x)$ in $\zeta^i$. It is sometimes called a discrete Fourier transform of $f$. In the following, we will prefer  the vectorial version of the Mattson-Solomon polynomial, which is 
$$\hat{f}:=(F_1,F_2,\ldots,F_n)=(f(\zeta),f(\zeta^2),\ldots,f(\zeta^n)).$$

\subsection{An invariant of fields}
We introduce here, following \cite{EKL}, the invariant of fields $$\mu(\F_q,n) := \min \{ d(C(f))+\dim C(f) \mid f \in R(\F_q,n), f\neq 0 \}.$$
By the Singleton bound, $\mu(\F_q,n)\leq n+1$ for any $n$. Moreover, equality holds if $n$ is prime and $q$ is a primitive root modulo $n$ or if $q$ is a power of $n$ (\cite[Propositions 4.3. and 4.4.]{EKL}). 

\begin{rem}
Note that $\mu(\F_q,n)=n+1$ if and only if all cyclic codes of length $n$ over $\F_q$ are MDS.
\end{rem}

\begin{rem}
Note that, as observed in \cite{EKL}, if we consider the complex field $\mathbb{C}$ instead of $\F_q$, then the  \textbf{uncertainty principle for simple cyclic group} $($proved for example in \cite{GGI,T}$)$ may be reformulated as follows: $\mu(\mathbb{C},p)=p+1$ for any prime $p$.
\end{rem}

In next sections we aim to investigate analogues of the uncertainty principle over finite fields.

\section{Strong version of UP}

The following version of UP is the one stated in \cite{EKL}. 

{\defi A finite field $\F_q$ satisfies the {\bf (strong) uncertainty principle} if, for all primes $p$, we have $\mu(\F_q,p)= p+1.$}

\medskip

As we have already mentioned, in \cite[Prop. 4.3]{EKL} it is shown that $\mu(\F_q,p)=p+1$ if $q$ is primitive mod $p$ (under this hypothesis, there exists exactly $3$ cyclic codes of length $p$ over $\F_q$ and of positive dimension, that are all trivial MDS codes). We show below that is almost the only case.

{\thm \label{MDS} Assume the MDS conjecture {\rm \cite[Res. Prob. 11.4]{MS}}.
	If $q$ is not primitive modulo $p$ and if $p>q+2$ then $\mu(\F_q,p)\leq p.$}

\begin{proof}
	By the hypothesis we know that there are polynomials $f | x^p-1$ in $R(\F_q,p)$ such that $1<\deg(f)<p-1.$ 
	Let $[p,k,d]$ be the parameters of the cyclic code $C(f).$
	Let $g\in C(f)$ of weight $d.$ The code $C(g)\subseteq C(f)$ is certainly not the repetition code, since $k>1.$
	Its parameters $[p,k'\le k,d]$ satisfy 
	$d+k'\ge \mu(\F_q,p)$ and if, arguing by contradiction,  $\mu(\F_q,p) > p,$ we see that
	$d \ge p-k'+1,$ entailing that
	$C(g)$ is MDS. But we know, by \cite[Chapt. 11]{MS}, that MDS codes of parameters $[N,K,D]$ with $1<K<N-1$ only exist for lengths at most $q+2.$
	This is the so-called MDS Conjecture that is now proved in many cases \cite{B,B+}.
	Note that codes of parameters $[N,1,N]$ and $[N,N-1,2]$ exist for all lengths $N.$
\end{proof}

{\rem \label{remMDS} A similar $($slightly weaker$)$ result holds unconditionally, since it is well-known that nontrivial MDS codes have length at most $2q-2$ $($see for example \cite[Corollary 7.4.4]{HP}$)$. So, with the same arguments we can prove that $\mu(\F_q,p)\leq p$ if $q$ is not primitive modulo $p$ and if $p>2q-2$.}

{\coro No finite field satisfies the $($strong$)$ uncertainty principle.}

\begin{proof}
    Suppose that $\F_q$ satisfies the $($strong$)$ uncertainty principle.
    Then Theorem \ref{MDS} would imply that all $p>q+2$ (or eventually $>2q-2$, if we refer to Remark \ref{remMDS}) are necessarily such that $q$ is primitive modulo $p$. But we know that this is not possible: it is enough to consider all primes $p$ such that $q$ is a quadratic residue (so that $q$ cannot be primitive modulo $p$) and observe that, by quadratic reciprocity, these correspond to $p$ being in some non-empty set of residue classes modulo $4q$ (so that they are infinitely many by Dirichlet's theorem).
\end{proof}

\section{Weak version of UP}
The following is Definition 5.3 in \cite{EKL}.

{\defi Let $0<\epsilon<\lambda\leq 1$. A finite field $\F_q$ satisfies the {\bf $(\epsilon,\lambda)$-uncertainty principle}, if 
	\begin{eqnarray}
	\mu(\F_q,p)&> &\lambda p\\
	ord_p(q) &<& \epsilon p,
	\end{eqnarray}
for infinitely many primes $p$, where $ord_p(q)$ is the order of $q$ dans $\F_p^\ast$.
}

\medskip

In \cite[Th. 5.4]{EKL} it is shown that finite fields satisfying this definition enjoy sequences of asymptotically good cyclic codes. Intuitively, (1) guarantees to get codes with a large minimum distance, whereas (2) guarantees to get codes with a large dimension.

In the following table, we show some values of $\mu(\F_2,p)$, for small primes $p$, omitting those for which $2$ is primitive modulo $p$.
$$\begin{array}{|c|c|c|c|c|c|c|c|c|c|c|c|c|c|c|c|} \hline
p&7 & 17 & 23& 31  & 41& 43 & 47 & 71& 73 & 79& 89& 97\\ \hline
\mu(\F_2,p)& 7 & 14  & 19 & 20 & 30& 28& 35& 47& 37 & 55& 45 & 64\\ \hline
\end{array}$$


In the following proposition, we get a restriction on possible values of $\lambda$ for finite fields satisfying the principle above.

{\prop If $\F_q$ satisfies the $(\epsilon,\lambda)$-uncertainty principle then $\lambda<\frac{q-1}{q}.$ 
}

\begin{proof}
By combining \cite[Th. 5.4]{EKL} with the same argument as in the proof of Theorem \ref{MDS}, we see that under the hypothesis, there 
are sequences of cyclic codes of length $p $ over $\F_q,$ of rate $R$ and relative distance $\delta$ such that
$$p \lambda < \mu(\F_q,p) < p\delta +pR.$$ In particular this implies that 

$$\lambda < \min \{ \delta+ \alpha_q(\delta) \mid \delta \in (0,1)\} , $$
where $\alpha_q(\delta)$ is the largest possible rate of a code of relative distance $\delta.$
But we know, by the asymptotic Plotkin bound \cite[Th.2.10.2]{HP}, that 
\begin{itemize}
	\item for $0<\delta < \frac{q-1}{q},$ we have $\alpha_q(\delta)< 1 -\frac{q\delta}{q-1},$ and
	\item for $\frac{q-1}{q} \le \delta <1,$ we have $\alpha_q(\delta)=0.$
\end{itemize}
It follows that the function $\delta+ \alpha_q(\delta)$ is $\le f(\delta)$ where 
\begin{itemize}
	\item  $f(\delta)= 1 -\frac{\delta}{q-1},$ for $0<\delta < \frac{q-1}{q},$  and
	\item $ f(\delta)=\delta  $ for $\frac{q-1}{q} \le \delta <1.$ 
\end{itemize}
Thus the minimum of $f(\delta)$ for $\delta \in (0,1)$ is met at $\delta=\frac{q-1}{q},$ and equals $\frac{q-1}{q}.$
\end{proof}

\section{Naive version of UP}

We prove here a finite field version of a result due to Donoho et al. \cite{D} in characteristic zero. Throughout this section we assume $(n,q)=1$ and we let $\zeta$ denote a primitive root of
unity of order $n$ in the algebraic closure of $\F_q$. 

{\prop[Naive UP] \label{donoho} For any $f\neq 0$ in $\F_q^n,$   
	$$ w_H(f)\cdot w_H(\hat{f})\ge n,$$
where $\hat{f}:=(f(\zeta), f(\zeta^2),\dots, f(\zeta^n))$ is the vectorial version of the Mattson-Solomon polynomial.
}

\begin{proof}
	Let $w:=w_H(f)$. By the BCH bound, $\hat{f}$ cannot have $w$ consecutive zeros.
	
	Suppose first that $w$ divides $n.$
	 Partition the set $\{1,\dots,n\}$ into
	$ n/w$ intervals of consecutive indices of length $w.$ In each of these intervals there is at least one index where $\hat{f}$ is nonzero. 
	Thus, we have exhibited $ n/w$ nonzeros of $\hat{f}.$ The desired inequality follows in that case.
	
	Equality holds if $\hat{f}$ has exactly one nonzero for each interval. Moreover, these nonzeros must be equally spaced, since otherwise there would be more than $w$ consecutive zeros between some pairs of nonzero elements of $\hat{f}$. 
	
	If $w$ does not divide $n$, then there is no way of distributing fewer than $\lceil n/w\rceil$ nonzero elements among $n$ places without leaving a gap of at least $w$ consecutive zeros. Thus $w_H(\hat{f})\geq \lceil n/w\rceil$.
\end{proof}

\noindent {\bf Remark:} The constant $n$ is best possible in view of the example of $f$ equal to the all-one vector (in this case $w_H(\hat{f})=1).$ Note also that a sharper bound has been very recently proved in \cite{FHX} by using van Lint-Wilson
bound \cite{LW}. \\

This can be reformulated in terms of cyclic codes as follows: for any $f \in R(\F_q,n)$, $f\neq 0$, 
$$d(C(f))\cdot \dim C(f) \geq n.$$

This allows to prove the following result, whose proof is technical and relegated to an appendix.

{\thm \label{2QC} For every real number $0<\alpha <1/2$  there are sequences of cyclic codes of rate $R$
	with minimum distance $\Omega(n^{\alpha}).$}
	
\medskip

\medskip

\noindent \textbf{Remark:} for $R\leq 1/2$, the square root bound on the minimum distance of quadratic residue codes (see for example \cite[Chap. 16, Th. 1]{MS}) gives an explicit construction of cyclic codes with asymptotic minimum distance bounded below by the square root of the length. However, for $R>1/2$ our result is the best, to our knowledge.

\section{Connection with Ramsey Theory}
In \cite{QRS} a connection between the uncertainty problem over finite fields and Ramsey Theory is pointed out. We give here a slight generalization, and a 
reinterpretation in terms of Coding Theory. We require a pair of definitions.

{\defi An {\bf arithmetic progression} of length $m$ in $\Z/n\Z$ is any subset of the form $\{a+kb \mid k\in\{0, \dots, m-1\}\}$ with $b \neq 0.$} 

{\defi The {\bf Szemer\'edi function} $r_m(n)$ is the largest size of a subset of $\Z/n\Z$ not containing an arithmetic progression of length $m.$}

{\prop \label{Ram} For $p$ prime such that $(q,p)=1$, we have \begin{center}
    $\mu(\F_q,p)\ge \min \{ m+p-r_m(p) \mid 1\le m \le p\}.$
\end{center}}
\begin{proof} Let $f \in R(\F_q,p).$ If $f$ has weight $m$, then, by BCH bound again, among the zeros of $f$ there cannot be $m$ consecutive powers of an $p$-th primitive root of unity. So $\{i\mid f(\zeta^i)=0, \zeta^p=1, \zeta\neq 1\}$ is a subset of $\Z/p\Z$ not containing an arithmetic progression of length $m$. Hence the number of zeros of $f$ is bounded above by $r_m(p)$. Then 
$$\begin{array}{rcl}\mu(\F_q,p) &=& \min \{ w_H(f)+w_H(\hat{f}) \mid  f \in R(\F_q,p),f\neq 0 \}\\& \geq &  \min \{w_H(f)+p-r_{w_H(f)}(p) \mid f \in R(\F_q,p),f\neq 0 \} \\ & \geq &  \min \{m+p-r_m(p) \mid 1\leq m\leq p \}.\end{array}$$
\end{proof}

\noindent {\bf Remark:} if $p$ is not prime, Proposition \ref{Ram} is not true. For example, $\mu(\F_2,9)=6$, but $\min \{ m+9-r_m(9) \mid 1\le m \le 9\}=8$. This is due to the fact that powers of $9$-th primitive root of unity may be $3$-rd root of unity.\\

\noindent {\bf Remark:} to our best knowledge, the function $r_m(p)$ is only known for $m$ fixed and $ p \to \infty$ \cite{TV}.
The fact that then $r_m(p)=o(p)$ is the celebrated Szemer\'edi Theorem. Any result on $r_m(p)$ when $m$ grows proportionally to $p$ with $p\to \infty$ would impact on the UP over finite fields.\\

We can generalize further by replacing arithmetic progression by their 2D analogues that is to say sets of the shape 
$$A(\delta,s)=\{a+kb+rc \mid k\in \{0, \dots, \delta-2\} \text{ and } \, r\in\{0, \dots, s\}\},$$ 
for $b,c$ coprime with $n$. Define then the function $r_{\delta,s}(n)$ as the largest 
size of a subset of $\Z/n\Z$ not containing an $A(\delta,s).$

{\prop For $p$ prime such that $(q,p)=1$, we have $$\mu(\F_q,p)\geq \min \{\delta+s-1+p-r_{\delta,s}(p) \mid \delta\in \{2,\ldots,p\} \text{ and } s\in\{0,\ldots,p-\delta\}\}.$$
}
\vspace{-5mm}
\begin{proof}
	The proof is the same as that of Proposition \ref{Ram}, up to the replacement of the BCH bound by the Hartmann-Tzeng bound \cite[p. 206]{MS}.
\end{proof}

\section{Conclusion and Open Problems}
In reaction to the recent papers \cite{EKL} and \cite{SWS}, we have considered the uncertainty principle when the Fourier transform takes its values over
finite fields. Exploring the connection with MDS codes, we prove that no finite field satisfies the strong version of UP introduced in \cite{EKL}. The weak version remains conjectural and we prove that 
it can only hold if $\lambda <\frac{q-1}{q}$. This should not discourage the researchers to try and prove
the weak version of UP for some values of $\lambda$ respecting this bound.

The analogue of the DFT UP of Donoho-Stark \cite{D}, which we called naive UP, allowed us to construct long cyclic codes of length $n$ and minimum distance 
$\Omega(n ^\alpha),$ where $0<\alpha<1/2.$ The proof is technical and relegated to an appendix. More suggestively, the arguments used to prove the naive UP yields an alternative proof 
of the results of \cite{SWS}, based on the BCH bound on the minimum distance of cyclic codes. A generalization based on the 
Hartmann-Tzeng bound has been sketched out. 

More generally, it would be worthy to generalize all these results to abelian codes.

\section*{Acknowledgement} The authors are grateful to Alexis Bonnecaze for programming help and to Pieter Moree for the fruitful discussions on the topic.


\medskip

\section*{Appendix: Proof of Theorem \ref{2QC}}
We construct a sequence of $q$-ary cyclic codes of rate $0<R<1,$ and designed minimum distance.
Let $p$ be an arbitrary prime, and write $n~=~q^p~-~1.$ If $x$ is an indeterminate then, from finite field theory \cite[Chap. 4, Th. 10]{MS}, we know that
$$x^n-1=\prod_{a\in \F_q^\ast}(x-a)\cdot \prod_{i=1}^sf_i,$$ where $f_i$ runs over all irreducible polynomials in $x$ of degree $p,$ and where $n=q-1+sp.$
Let $g_I= \prod_{i\in I}f_i$ with $|I|=s'=\lfloor s(1-R)\rfloor.$ Then the dimension of the cyclic code $C(g_I)$ of generator $g$ is $n-ps',$ and it can be checked that
$$\frac{n-s'p}{n} \to R$$ when $ p\to \infty.$

We need to control the intersections of the $C(g_I)$'s when $I$ varies.

{\lem \label{mu} Let $r\neq 0$ be an arbitrary vector in $\F_q^n,$ of  Hamming weight $\le n^\alpha$ for some $0<\alpha<1$. There at most
	$\Lambda_n(1+o(1))$ polynomials $g_I$ with $|I|=s'$
	such that $r\in C(g_I),$ where $\Lambda_n=2^{\left(\frac{\delta n-n^{1-\alpha}}{\delta p}\right)H(R)}.$}

\begin{proof}
	
	The number $N$ of indexes $I$ such that $r\in C(g_I),$ equals the number of $s'$-sets of $f_i$'s which divide $r(x)$ in polynomial notation.
	
	Let $\zeta$ be a primitive $n$-th root of unity in the algebraic closure $\mathbb{K}$ of $\F_q$, and $\widehat{r}=(r(\zeta), r(\zeta^2),\dots,r( \zeta^n))$. We have
	$$N\leq {{\lfloor \frac{Z(r)}{p} \rfloor }\choose s'},$$
	where
	$$ Z(r):=| \{ \omega\in \mathbb{K} \mid \omega^n=1 \text{ and }  r(\omega)=0\}|=n-w_H(\widehat{r}).$$
	By Lemma \ref{donoho}, $w_H(\widehat{r})\ge n^{1-\alpha},$ so that $Z(r)\le n-n^{1-\alpha}.$ Thus
	
	$$ N\leq {{\left\lfloor \frac{n-n^{1-\alpha}}{p} \right\rfloor }\choose s'}\le 2^{\left(\frac{ n-n^{1-\alpha}}{p}\right)H(R)},$$
	where the upper bound is a consequence of \cite[Chap. 10, Lemma 8]{MS}.
\end{proof}

\begin{proof}[Proof of Theorem \ref{2QC}]
	The number of possible $g_I$'s is $${ s \choose s'}\sim \frac{2^{sH(R)}}{\sqrt{2\pi s R(1-R)}}$$ for $s \to \infty,$ by
	Stirling's approximation of the factorial.
	
	If this number is greater than the product of $\Lambda_n$ by the volume of the Hamming ball of radius 
	$n^\alpha$ in length $n,$
	then there are codes $C(g_I)$ with minimum distance at least $n^\alpha.$ 
	
	The volume of the
	Hamming ball of radius $n^\alpha$ is bounded above by 
	$$(1+\lfloor n^\alpha\rfloor){n \choose \lfloor  n^\alpha\rfloor}(q-1)^{\lfloor n^\alpha\rfloor}$$
	(see the proof of \cite[Lemma  2.10.3]{HP}), which is bounded above by a quantity asymptotically equivalent to
	$$\frac{2^{nH(n^{\alpha-1})+ \log_2(n^\alpha)+ n^\alpha\log_2(q-1)}}{\sqrt{2\pi  n^{\alpha}(1- n^{\alpha-1})}}$$
	$$\sim \frac{1}{\sqrt{2\pi}}\cdot n^{-n^\alpha(\alpha-1)+\frac{\alpha}{2}}\cdot e^{n^\alpha-n^{2\alpha-1}}\cdot (q-1)^{n^\alpha}$$
	for $n \to \infty$. So, applying Lemma \ref{mu}, the mentioned inequality happens if
	
	$$2^{\left(\frac{-(q^p-1)^{1-\alpha}+(q-1)}{p}\right)H(R)}\cdot   (q^p-1)^{-(q^p-1)^\alpha(\alpha-1)+\frac{\alpha}{2}}\cdot $$
	$$\cdot
	e^{(q^p-1)^\alpha-(q^p-1)^{2\alpha-1}}\cdot (q-1)^{(q^p-1)^\alpha} \cdot \left(\frac{q^p-q}{p}\right)^{1/2} \leq \frac{1}{\sqrt{R(1-R)}}.$$
	We can write the last inequality as $e^{f_{\alpha,q,R}(p)}\leq \frac{1}{\sqrt{R(1-R)}}$, with
	$$f_{\alpha,q,R}(p)=(1-\alpha)\ln(q^p-1)(q^p-1)^\alpha-\ln(2)H(R)\frac{q^p}{p(q^p-1)^{\alpha}}+o(pq^{\alpha p})$$
	for $p\to \infty$, so that $f_{\alpha,q,R}(p)\to -\infty$ for $p\to \infty$ if $\alpha<1/2$, and it grows to $\infty$ otherwise.
\end{proof}
\end{document}